\documentclass[reqno, 11pt]{amsart}
\usepackage{mathtools}
\usepackage{amsmath}
\usepackage{amssymb}
\usepackage{yhmath}
\usepackage{graphicx}
\usepackage{ mathrsfs }
\usepackage{bbm}
\usepackage{xcolor}
\usepackage{tikz-cd}
\usepackage{tikz}
\usetikzlibrary{patterns}
\usepackage{hyperref}

\DeclareMathAlphabet{\mathpzc}{OT1}{pzc}{m}{it}

\usepackage{thmtools}
\usepackage{thm-restate}

\usepackage{caption}

\newtheorem{theorem}{Theorem}[section]

\newtheorem*{claim*}{Claim}

\newtheorem{lemma}[theorem]{Lemma}

\newtheorem{cor}[theorem]{Corollary}

\newtheorem{prop}[theorem]{Proposition}%https://www.overleaf.com/project/616d911a46ea936b949fdb96
\newtheorem{Thm}[theorem]{Theorem}
\newtheorem{Con}[theorem]{Conjecture}

\theoremstyle{definition}
\newtheorem{Def}[theorem]{Definition}

\newtheorem{rmk}[theorem]{Remark}
\newtheorem{Rmk}[theorem]{Remark}

\numberwithin{equation}{section}

\newcommand{\op}{\operatorname}

\newcommand{\bb}{\mathbb}

\newcommand{\be}{\begin{equation}}
\newcommand{\ee}{\end{equation}}
\newcommand{\Ga}{\Gamma}
\renewcommand{\c}{\mathbb C}
\newcommand{\R}{\mathbb R}

\newcommand{\ga}{\gamma}

\newcommand{\La}{\Lambda}
\newcommand{\inte}{\op{int}}
\newcommand{\ba}{\backslash}

\newcommand{\cal}{\mathcal}
\newcommand{\br}{\mathbb R}
\newcommand{\SO}{\op{SO}}

\newcommand{\PSL}{\op{PSL}}
\newcommand{\F}{\cal F}

\newcommand{\bH}{\mathbb H}

\newcommand{\G}{\Gamma}

\renewcommand{\frak}{\mathfrak}

\newcommand{\e}{\varepsilon}

\renewcommand{\L}{\mathcal L}
\newcommand{\fa}{\mathfrak a}

\renewcommand{\S}{\mathbb S}

\newcommand{\so}{\SO^\circ(n,1)}

\newcommand{\s}{\sigma}

\newcommand{\fg}{\frak g}

\begin{document}

\title[Tent property and applications]{Tent property of the growth indicator functions
and applications}

\author{Dongryul M. Kim}
\address{Department of Mathematics, Yale University, New Haven, CT 06511}
\email{dongryul.kim@yale.edu}

\author{Yair N. Minsky}
\address{Department of Mathematics, Yale University, New Haven, CT 06511}\email{yair.minsky@yale.edu}

\author{Hee Oh}
\address{Department of Mathematics, Yale University, New Haven, CT 06511}
\email{hee.oh@yale.edu}
\thanks{Minsky and Oh are partially supported by the NSF}

\begin{abstract} 
Let $\Gamma$ be a Zariski dense discrete subgroup of a connected semisimple real algebraic group $G$. Let $k=\operatorname{rank} G$. Let $\psi_\Ga:\mathfrak{a} \to \mathbb{R}\cup \{-\infty\}$ be the growth indicator function of $\Ga$, first introduced by Quint.
In this paper, we obtain the following pointwise bound of $\psi_\Gamma$: for all $v\in \mathfrak{a}$,
$$ \psi_\Gamma(v) \le \min_{1\le i\le k} \delta_{\alpha_i} \alpha_i(v) $$
where $\Delta=\{\alpha_1, \cdots, \alpha_k\}$ is the set of all simple roots of $(\mathfrak{g},\mathfrak{a})$ and
$0<\delta_{\alpha_i}\le \infty$ is the critical exponent of $\Gamma$ associated to $\alpha_i$.
 When $\Gamma$ is $\Delta$-Anosov, there are precisely $k$-number of directions where
the equality  is achieved, and the following strict inequality holds for $k\ge 2$: for all $v\in \mathfrak{a}-\{0\}$,
 $$\psi_\Gamma(v) <\frac{1}{k}\sum_{i=1}^k \delta_{\alpha_i} \alpha_i (v).$$ We discuss applications  for 
 self-joinings of convex cocompact subgroups in $\prod_{i=1}^k \operatorname{SO}(n_i,1)$ and 
 Hitchin subgroups of $\operatorname{PSL}(d, \mathbb{R})$. In particular, 
 for a Zariski dense Hitchin subgroup $\Gamma<\text{PSL}(d, \mathbb{R})$,
 we obtain that
 for any $ v=\operatorname{diag}(t_1, \cdots, t_d)\in \mathfrak{a}^+$,
$$\psi_\Gamma (v) \le \min_{1\le i\le d-1} (t_i -t_{i+1}). $$

\end{abstract}

\maketitle
\tableofcontents

\section{Introduction}
Let $G$ be a connected semisimple real algebraic group.
We let  $P=MAN$ be a minimal parabolic subgroup of $G$ with a fixed Langlands decomposition, where  $A$ is a maximal real split torus of $G$, $M$
is the maximal compact subgroup centralizing
$A$ and $N$ is the unipotent radical of $P$. Let $\fg=\op{Lie} G$, $\fa=\op{Lie} A$ and $\fa^+$ denote the positive Weyl chamber so that $\log N$ consists of positive root subspaces. 
Let $K$ be a maximal compact subgroup so that the Cartan decomposition $G=K(\exp \fa^+) K $ holds.
Let $\mu:G\to \fa^+$ denote the Cartan projection map defined by the condition
$\exp \mu(g)\in KgK$ for all $g\in G$.
Let $\Ga<G$ be a Zariski dense discrete subgroup. We denote by $\L\subset \fa^+$
the limit cone of $\Ga$, which is the asymptotic cone of $\mu(\Ga)$. It is a convex cone with non-empty interior \cite{Ben}.

Following Quint \cite{Quint1}, the growth indicator function $\psi_\Ga : \fa \to\bb R\cup\{-\infty\}$  is defined as follows: choose any  norm $\|\cdot\|$ on $\fa$.
For an open cone $\cal C$ in $\fa$, let $\tau_{\cal C}$ denote the abscissa of convergence of  $\sum_{\ga\in\Ga,\,\mu(\ga)\in\cal C}e^{-s \|\mu(\ga) \|}$ (that is, the infimum of the set of $s$ for which the series converges). Now 
for any non-zero $v\in \fa$, let 
\begin{equation}\label{grow3}\psi_{\Gamma}(v):=\|v\| 
\inf_{v\in\cal C} \tau_{\cal C}\end{equation}
where the infimum is over all open cones $\cal C$ containing $v$, and let $\psi_\Gamma(0)=0$.  The definition of $\psi_\Ga$ does not depend on the choice of a 
norm on $\fa$.
Note that $\psi_\Gamma =-\infty$ outside $\L$. 
Quint showed that $\psi_\Ga$ is a concave upper-semi continuous function satisfying $\L=\{\psi_\Ga\ge 0\}$ and $\psi_\Ga>0$ on the interior
$\inte \L$.

The main aim of this paper is to present a pointwise bound for the growth indicator function together with some applications. 
Throughout the paper, for any non-negative function $f$ on $\fa^+$,
we denote by $$0\le \delta_{\Ga, f}\le \infty$$ or simply, $\delta_{f}$,
the critical exponent of $\G$ with respect to $f$, that is, the abscissa of convergence of the series $\sum_{\ga \in \Ga} e^{-s f(\mu(\ga))}$.

Let 
$$\Delta=\{\alpha_1, \cdots, \alpha_k\}$$ denote the set of simple roots for $(\frak g, \frak a^+)$.

 \begin{Def}[Tent function]\label{tttdef} Let $\Ga<G$ be a Zariski dense discrete subgroup
 with $\delta_{\Ga, \alpha_i}<\infty$ for some $1\le i\le k$.
We define a tent function $T_\Ga:\fa\to [0, \infty)$ by
$$T_\Ga (v):= \min_{1\le i\le k} \delta_{\Ga, \alpha_i}\cdot  \alpha_i(v)  .$$
\end{Def}

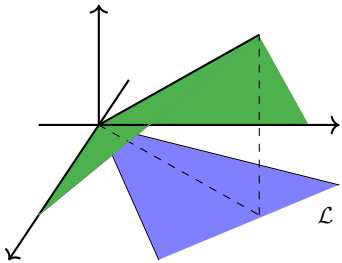
\begin{figure}[ht]
\centering
\begin{tikzpicture}[scale=0.8, every node/.style={scale=0.8}]
\filldraw[draw=white, fill = blue!50, opacity=0.2] (0, 0) -- (4, -1) -- (1, -2.25) -- (0, 0);
\draw (0, 0) -- (4, -1);
\draw (0, 0) -- (1, -2.25);
%\draw[->] (0, 0) -- (1.5, -0.6);
%\draw[->] (0, 0) -- (1, -1.2);
%\draw (1.5, -0.6) node[right] {$u_{e_1}$};
%\draw (1, -1.2) node[right] {$u_{e_2}$};
\draw (3.5, -1.5) node[right] {$\L$};
\filldraw[draw = white, fill = green!40!gray, opacity=0.3] (0, 0) -- (3.5, 0) -- (8/3, 1.5) -- (0, 0);
\filldraw[draw=white, fill= green!40!gray, opacity = 0.3] (8/3, 1.5) -- (0, 0) -- (-1, -1.5);
\draw[->, thick] (0, 0) -- (0, 2);
\draw[->, thick] (-1, 0) -- (4, 0);
\draw[->, thick] (1/2, 3/4) -- (-1.5, -2.25);
\draw[dashed] (0, 0) -- (8/3, -1.5);
\draw[dashed] (8/3, -1.5) -- (8/3, 1.5);
\draw[draw=green!30!gray] (8/3, 1.5) -- (-1, -1.5);
\draw[draw=green!30!gray] (8/3, 1.5) -- (-1, -1.5); (8/3, 1.5) -- (3.5, 0);
\draw[thick] (8/3, 1.5) -- (0, 0);
\end{tikzpicture}
\caption{\small{Tent on the limit cone}} \label{fig.tent}
\end{figure}

We obtain the following tent property of the growth indicator function: 
\begin{Thm}[Tent property]\label{m}
For any Zariski dense discrete subgroup $\Ga<G$ such that $\min_{1\le i\le k} \delta_{\Ga, \alpha_i}<\infty$, we have
$$ \psi_\Ga (v)\le T_\Ga(v)\quad\text{for all $v\in \fa$}. $$
Moreover, when $\delta_{\Ga, \alpha_i}<\infty$, there exists $v_i \in \L - \{0\}$
such that $\psi_\Ga(v_i)=T_\Ga(v_i)=\delta_{\Ga, \alpha_i}\alpha_i(v_i)$.
\end{Thm}
\begin{rmk}\label{Ka}
\begin{enumerate}
    \item Denote by $\pi_G$ the half-sum of all positive roots of $(\fg, \fa^+)$ counted with multiplicity. Then for any discrete subgroup $\Ga<G$,
    we have $\psi_\Ga\le 2\pi_G$  \cite[Thm. IV.2.2]{Quint1}.
    
    \item  If $G$ has property $\op{(T)}$ and $\Ga$ is of infinite co-volume, then
    $\psi_\Ga \le 2\pi_G -\Theta$ where $\Theta$ is the half-sum of a maximal strongly orthogonal system (\cite{Q3}, \cite{Oh}, see also \cite[Thm. 7.1]{LO1}). Our bound  in Theorem \ref{m} provides a sharper bound for Hitchin subgroups; see Remark \ref{comparePS}.
\end{enumerate}
\end{rmk}

For a non-empty subset $\theta\subset \Delta$, 
a finitely generated  subgroup $\Ga<G$ is called a $\theta$-Anosov subgroup
if  there exist constants $C,C'>0$ such that for all $\ga\in\Ga$ and all $\alpha_i\in \theta$,
\be \label{eqn.Anosovdef}
\alpha_i(\mu(\ga))\ge C|\ga|-C'
\ee
 where $|\ga|$ denotes the word length of $\ga$ with respect to a fixed finite symmetric set of generators of $\Ga$.
 The notion of Anosov subgroups was first introduced by Labourie for surface groups \cite{La}, and was extended to general word hyperbolic groups by Guichard-Wienhard \cite{GW}.
Several equivalent characterizations have been  established, one of which is the above definition (see \cite{GGKW}, \cite{KL}, \cite{KLP1}, \cite{KLP2}). Anosov subgroups are regarded as natural generalizations of convex cocompact subgroups of rank one groups.

For a $\theta$-Anosov subgroup $\Ga<G$, it follows from \eqref{eqn.Anosovdef} that for some constant $C>0$,
$$\max_{\alpha_i\in \theta} \delta_{\Ga,\alpha_i}\le C \log \#\mathsf S <\infty$$   where $\mathsf S$ is a fixed finite generating set of $\Ga$. Therefore
 Theorem \ref{m} applies to any Zariski dense subgroup contained in some $\theta$-Anosov subgroup of $G$.

For $\Delta$-Anosov subgroups, we obtain the following sharper result:
\begin{theorem} \label{thm.Anosov}\label{m4}
 Let $\Ga$ be a Zariski dense $\Delta$-Anosov subgroup of $ G$. The following hold:
\begin{enumerate}
    
\item For each $1\le i\le k$, there exists a unique  $v_i \in \inte \L$ such that 
$\alpha_i(v_i)=1$ and $\psi_\Ga (v_i)=\delta_{\Ga, \alpha_i}$. \label{(2)}

\item For $v\in \fa-\{0\}$, we have
$\psi_\Ga(v)\le T_\Ga(v)$ where equality holds
if and only if $v=cv_i$ for some $1\le i\le k$ and $c>0$. \label{(3)}

\item If $k=\op{rank} G \ge 2$, then  $$ \psi_{\Ga} < \frac{1}{k} \sum_{i=1}^k \delta_{\Ga,\alpha_i} \alpha_i.$$ \label{(4)}

\end{enumerate}
\end{theorem}

 When $\Gamma$  is $\Delta$-Anosov, $\psi_\Ga$ is strictly concave\footnote{Since $\psi_\Ga$ is homogeneous, the strict concavity of $\psi_\Ga$ is equivalent to saying that $\psi_\Ga (v+w) > \psi_\Ga(v)+\psi_\Ga(w)$ for all $v, w\in \inte \L$ in different directions} in $\inte \L$ by (\cite[Thm. A]{Samb3}, \cite[Prop. 4.11]{SP}).
Therefore by the convexity of the unit norm ball $\{\|v\|\le 1\}$, there exists a unique unit vector $u_{\Ga, \|\cdot\|}\in \fa^+$, called the direction of maximal growth, such that  $\psi_\Ga(u_{\Ga, \|\cdot\|} )=\max_{\|v\|=1}\psi_{\Ga}(v) .$
By \cite[Coro. III.1.4]{Quint1}, we have
 \be\label{dddd}
\delta_{\Ga, \|\cdot\|} =\psi_\Ga(u_{\Ga, \|\cdot\|} ).
\ee

\begin{cor} \label{innerprod} Let $k=\op{rank} G \ge 2$. Let $\Ga$ be a Zariski dense $\Delta$-Anosov subgroup of $ G$.
For any  norm $\|\cdot\|$ on $\fa$ induced from an inner product which is non-negative on $\fa^+$, we have
$$\delta_{\Ga, \|\cdot\|} < \min_{1\le i\le k} \delta_{\Ga, \alpha_i} \cdot \alpha_i(u_{\Ga, \|\cdot\|}).$$ 
\end{cor}

In view of the above discussion,  any upper bound on $\delta_{\Ga, \alpha_i}$ for any $\alpha_i\in \Delta$ provides an explicit pointwise upper bound on $\psi_\Ga$. We discuss some examples of $\Delta$-Anosov subgroups.

\subsection*{Self-joinings of hyperbolic manifolds}
For $1\le i\le k$, consider the hyperbolic space $(\bH^{n_i},d_i)$, $n_i\ge 2$, with constant sectional curvature $-1$, and let $G_i=\SO^{\circ}(n_i, 1)=\op{Isom}^+(\bH^{n_i})$.
Let $G=\prod_{i=1}^kG_i$. 
Denote by $\alpha_i$ the simple root of $\frak g_i = \op{Lie} G_i$. Then $\Delta=\{\alpha_1, \cdots,
\alpha_k\}$ is the set of simple roots of $\fg$.
Via
the map $v\mapsto (\alpha_1(v), \cdots, \alpha_k(v))$, we may identify
$\fa=\br^k$ and $\fa^+=\{(v_1, \cdots, v_k)\in \br^k: v_i\ge 0\text{ for all $i$}\}$.

Let $\Sigma$ be a countable group and $\rho_i:\Sigma\to G_i$ be a faithful 
convex cocompact representation with Zariski dense image for each $1\le i\le k$.
Setting $\rho=(\rho_1, \cdots, \rho_k)$,
the self-joining $\Ga_\rho$ is defined as
the following  subgroup of $G$:
\be\label{grho} \Ga_\rho=\left(\prod_{i=1}^k \rho_i\right)(\Sigma)=\{ (\rho_1(\s), \cdots, \rho_k(\s))\in G: \s\in \Sigma\}.\ee 
We also assume that no two of $\rho_i$'s are conjugate, so that $\Ga_{\rho}$ is a Zariski dense discrete subgroup of $G$. The hypothesis on $\rho_i$'s implies that $\Ga_\rho$ is a $\Delta$-Anosov subgroup of $G$ (cf. \cite[Thm. 5.15]{GW}).

Fix $o_i\in \bH^{n_i}$. For each $1\le i\le k$, denote by $0<\delta_{\rho_i}\le \infty$ the critical exponent of $\rho_i(\Sigma)$, that is, the abscissa of convergence of the series $\sum_{\sigma \in \Sigma} e^{-s d_i(\rho_i(\sigma)o_i, o_i)}$. 
We also denote by  $\La_{\rho_i}\subset 
\S^{n_i - 1}$ the limit set of $\rho_i(\Sigma)$, which is the set of accumulation points of $\rho_i(\Sigma)o_i$ in the compactification
$\bH^{n_i}\cup \S^{n_i-1}$. These two notions are independent of the choice of $o_i\in \bH^{n_i}$. By Patterson \cite{Pa} and Sullivan \cite{Su},  we have
\be \label{ci}
\delta_{\rho_i} = \dim \La_{\rho_i}
\ee
where $\dim \La_{\rho_i}$ is the Hausdorff dimension of $\La_{\rho_i}$ with respect to the spherical metric $d_{\S^{n_i- 1}}$. We deduce from Theorem \ref{thm.Anosov}:
\begin{cor}\label{pregap} Let $\Ga_\rho<G$
be a Zariski dense subgroup
of $G=\prod_{i=1}^k \SO^{\circ}(n_i, 1)$, $n_i\ge 2$, as defined in \eqref{grho}. Assume $k\ge 2$. For any $v = (v_1, \cdots, v_k) \in \R^k$, we have $$ \psi_{\Ga_{\rho}}(v) < \frac{1}{k} \sum_{i = 1}^k \dim \La_{\rho_i} \cdot v_i.$$ In particular, 
we have $$\delta_{\Ga_{\rho},\|\cdot\|_{\op{Euc}}} < \frac{1}{k}{\left({\sum_{i = 1}^k (\dim \La_{\rho_i})^2}\right)^{1/2}} $$
where $\|\cdot\|_{\op{Euc}}$ denotes
 the standard Euclidean norm on $\br^k$. 
\end{cor}

 Let $\F=\prod_{i=1}^k\S^{n_i-1}$, which is the Furstenberg boundary of $G$.
 The limit set of $\Ga_\rho$ is the set of all accumulation points of an orbit $\Ga_\rho 
 (o_1, \cdots, o_k)$: 
\be\label{dll} \La_\rho=\left\{ (\xi_1, \cdots, \xi_k)\in \F: \begin{matrix}
\text{$\exists$ a sequence $\sigma_\ell\in \Delta$ s.t. $\forall$ $1 \le i \le k$,} \\
\xi_i=\lim_{\ell\to \infty} \rho_i(\sigma_\ell)(o_i)
\end{matrix} \right\} .\ee

In \cite{KMO}, we showed that 
\be\label{kmo}\dim \La_\rho=\max_{1\le i\le k} \dim \La_{\rho_i}\ee 
where the Hausdorff dimension of $\La_{\rho}$ is computed with respect to the Riemannian metric on $\F$ given by
$\sqrt{\sum_{1\le i\le k} {d_{\S^{n_i- 1}}}^2}$.  We deduce the following from Corollary \ref{pregap} and \eqref{kmo}:

\begin{cor}[Gap theorem]\label{rigid} For  $k\ge 2$, we  have
$$\delta_{\Ga_\rho, \|\cdot\|_{\op{Euc}}} < \frac{\dim \La_\rho}{\sqrt{k}}.
$$
 \end{cor}

The trivial bound for 
$\delta_{\Ga_\rho, \| \cdot \|_{\op{Euc}}} $ is given by
$\delta_{\Ga_\rho, \| \cdot \|_{\op{Euc}}} \le\min_i \delta_{\rho_i}  \le \dim \La_{\rho}$.
Hence Corollary \eqref{rigid} presents a strong gap for the value of $\delta_{\Ga, \|\cdot\|_{\op{Euc}}}$ from the trivial bound. 
This phenomenon is in contrast to the rank one case:
there exist convex cocompact (non-lattice) subgroups
$\Ga$ of $\SO^\circ(n,1)$ whose critical exponents $\delta_\Ga$ are arbitrarily close to $n-1$ (see e.g., \cite[Sec.6]{Mag} on the construction of McMullen). 

\begin{Rmk}  Let $\rho_1, \rho_2$ be two convex cocompact faithful representations into $\so=\op{Isom}^\circ (\bH^n)$ and $\rho=(\rho_1, \rho_2)$. 
 Note that $\Ga_\rho <\so\times \so$ is Zariski dense if and only if $\rho_1$ and $\rho_2$ are not conjugate by an element of $\op{Isom}^\circ (\bH^n)$. Hence Corollary \ref{rigid} can be interpreted as the following rigidity statement: we have
\be\label{gap3}\delta_{\Ga_\rho, \|\cdot\|_{\op{Euc}}}\le    \frac{n-1}{\sqrt 2}\ee 
and the equality holds
if and only if $\rho_1(\Sigma)$ and $\rho_2(\Sigma)$ are conjugate {\it lattices} of $\so$. This particular rigidity statement was recently extended  in
\cite{BC} even to geometrically finite representations.
 \end{Rmk}

In view of special interests in low dimensional hyperbolic manifolds which come with huge deformation spaces, we also formulate the following consequence 
of Corollary \ref{innerprod}, using the isomorphisms $\PSL(2, \c)\simeq \SO^\circ (3,1)$ and $\PSL(2, \br)\simeq \SO^\circ (2,1)$, the characterization of the critical exponent in \eqref{dddd}, and  the simple fact 
$ \sup\{\min(v_1, 2v_2): v_1^2+v_2^2=1\} = \frac{2}{\sqrt 5} $.

\begin{cor}\label{delta} Consider the metric on $\bH^2\times \bH^3$ given by
$d= \sqrt{d_{\bH^2}^2 + d_{\bH^3}^2} $.
For any non-elementary convex cocompact subgroup $\Gamma_0<\PSL(2, \br)$
and any non-elementary faithful convex cocompact Zariski dense representation
$\rho_0 : \Gamma_0 \to \PSL(2, \c)$,
the critical exponent of the group $\{(\gamma_0, \rho_0(\gamma_0)) : \gamma_0 \in \Ga_0 \}$
with respect to $d$ 
is strictly less than  $  \frac{2}{\sqrt 5}$.
\end{cor}

\subsection*{Hitchin representations.}
We discuss applications to Hitchin representations.
In $G=\PSL(d, \br)$, we have
$\fa^+=\{v=\text{diag}(t_1, \cdots, t_d): t_1\ge \cdots \ge t_d, \sum t_i=0\}$ and  $\alpha_i (v)= t_{i}-t_{i+1}$ for $1\le i\le d-1$.
Let $\Sigma$ be  a torsion-free uniform lattice of $\PSL(2,\br)$,
and $\pi_d$ denote the $d$-dimensional irreducible representation   $\PSL(2, \br)\to \PSL(d, \br)$, which is unique up to conjugation.
A Hitchin representation
$\rho:\Sigma\to \op{PSL}(d,\br)$ is a representation which
 belongs to the same connected component as $\pi_d|_{\Sigma}$ in the character variety $\op{Hom}(\Sigma, \PSL(d, \br))/\sim$ where the equivalence is given by conjugations. 
 
 We call the image of a Hitchin representation $\Ga:=\rho(\Sigma)$ a Hitchin subgroup of $G$.

A Hitchin subgroup is known to be a $\Delta$-Anosov subgroup of $\PSL(d, \br)$ by Labourie \cite{La}. 
 By the work of Potrie-Sambarino \cite[Thm. B]{SP}
(see also \cite[Coro. 9.4]{PSW2}), a Hitchin subgroup
$\Ga<\PSL(d, \br)$ satisfies:
\be\label{ps2}\delta_{\Ga, \alpha_i}=1 \quad\text{for all $1\le i\le d-1$}.\ee 
Together with this important result,
Theorems \ref{m} and \ref{thm.Anosov}  imply
the following:
\begin{cor} \label{hit} Let $d\ge 3$ and $\G<\op{PSL}(d,\br)$ be a Zariski dense Hitchin subgroup of $\op{PSL}(d,\br)$.
Then for any $ v=\op{diag}(t_1, \cdots, t_d)\in \fa^+$,
\begin{equation}\label{point1} \psi_\Ga (v) \le \min_{1\le i\le d-1} (t_i -t_{i+1}); \end{equation} 
\begin{equation} \psi_{\Ga}(v) < (t_1-t_d)/(d-1).\end{equation} 
\end{cor}

 This pointwise bound for $\psi_\Ga$ is sharper than
 the one from (\cite{Q3}, \cite{Oh}, \cite[Thm. 7.1]{LO1}), which for instance, for $d=3$,
 gives the upper bound $\frac{3}{2} (t_1-t_3)$ while the above corollary gives a bound $\frac{1}{2}(t_1-t_3)$. 

\begin{rmk} Following \cite{CZZ}, for any geometrically finite subgroup $\Sigma<\PSL(2, \br)$, a representation $\rho:\Sigma\to \op{PSL}(d,\br)$ is called cusped Hitchin if there exists a {\it positive} $\rho$-equivariant map from the limit set of $\Sigma$  to the space $\cal F$ of complete $d$-dimensional flags.   For a cusped Hithin subgroup  $\Gamma<\PSL(d, \br)$, i.e., the image of a cusped Hitchin representation of a geometrically finite $\Sigma<\PSL(2, \br)$, the inequality \begin{equation}
    \label{czz} \max_{1\le i\le d-1} \delta_{\Ga, \alpha_i}\le 1 \end{equation} was obtained, with equality only when $\Sigma$ is a lattice, by
Canary, Zhang and Zimmer \cite[Thm. 1.1]{CZZ}. Although $\Ga$ is not Anosov when $\Sigma$ is not convex cocompact,  Theorem \ref{m}, using \eqref{czz}, implies that
the pointwise bound \eqref{point1}
$\psi_\Ga (v) \le \min_{1\le i\le d-1} (t_i -t_{i+1})$, and hence $\psi_\Ga(v)\le (t_1-t_d)/(d-1) $, also holds for any Zariski dense cusped Hitchin subgroup $\Gamma$ of $\PSL(d,\br)$ as well.
 
\end{rmk} 
 
\begin{rmk}
The bound in Corollary \ref{hit} is stronger than \cite[Coro. 1.4]{SP} (also \cite[Thm. 1.1]{CZZ} for cusped Hitchin subgroups) in two aspects:
first, the bound for $\psi_\Ga$ given by \cite[Coro. 1.4]{SP} is weaker than 
$\frac{t_1-t_d}{d-1}$ and stated only for vectors inside a strictly smaller cone than the limit cone (see Remark \ref{comparePS} for details).
\end{rmk}

 \begin{rmk} 
The comparison of $\psi_\Ga$ with the half sum $\pi_G$ of positive roots is meaningful in view of Sullivan's theorem that for a convex cocompact subgroup $\Ga < \so$, the inequality $\delta_{\Ga} \le \pi_G=\frac{n - 1}{2}$ holds if and only if the bottom of the $L^2$-spectrum on $\Ga \ba \bH^n$ is given by $(n-1)^2/4$ and there exists no positive square-integrable harmonic function on $\Ga \ba \bH^n$ \cite[Thm. 2.21]{Su2}.
 
 Corollaries \ref{pregap} and \ref{hit}
 imply that $\psi_\Ga \le  \pi_G $ in their respective settings (even with the strict inequality).  In recent work \cite{EO}, these results were used to show that the quasi-regular representation $L^2(\Ga\ba G)$ is tempered and
 there exists no positive square-integrable harmonic function on the associated locally symmetric manifold.

 For any discrete subgroup 
$\Ga<G$, note that $\delta_{\Ga, \pi_G} \le 2$ as follows from Remark \ref{Ka}(1). 
We propose the following conjecture:
\begin{Con}\label{conj} Let $k=\op{rank} G\ge 2$.
If $\Ga$ is a $\Delta$-Anosov subgroup of $G$, then $$\delta_{\G,\pi_G} \le 1,$$ or equivalently $\psi_\Ga\le \pi_G$.
\end{Con} 
The equivalence is a consequence of \cite[Lemma III.1.3]{Quint1}.
 \end{rmk}

\subsection*{On the proofs}
 The proof of Theorem \ref{m} consists of two parts: first
 prove that each linear form
$\delta_{\Ga,\alpha_i}\alpha_i$ is tangent to $\psi_\Ga$ whenever $\delta_{\Ga, \alpha_i}<\infty$ and then  take the {\it minimum}! Although taking the minimum
seems a trivial step, the resulting tent function turns out to be quite useful, as discussed above.
The proof of Theorem \ref{thm.Anosov} is crucially based on special properties of $\psi_\Ga$ for $\Delta$-Anosov subgroups (see Theorem \ref{SPS}).

\subsection*{Organization} 
In section \ref{prelim}, we prove Theorem \ref{m}.
In section \ref{sec.dual}, we prove Theorem \ref{thm.Anosov}. 
In section \ref{sec.self}, we discuss applications of tent property of $\psi_{\Ga}$ to self-joining of hyperbolic manifolds. 

\subsection*{Acknowledgements}
We would like to thank Marc Burger and  Dick Canary  for useful comments 
and Andres Sambarino for pointing out
 some redundant rank restriction in our earlier version.

\section{Tent property} \label{prelim}
Let $G$ be a connected, semisimple real algebraic group of rank $k\ge 1$.
Let $\fg$ denote the Lie algebra of $G$, and decompose $\fg$ as $\mathfrak g=\mathfrak k\oplus\mathfrak{p}$, where $\frak k$ and $\frak p$ are the $+ 1$ and $-1$ eigenspaces of 
a fixed Cartan involution respectively. We denote by $K$ the maximal compact subgroup of $G$ with Lie algebra $\frak k$. We also choose a maximal abelian subalgebra $\fa$ of $\mathfrak p$. Let $A:=\exp \mathfrak a$. Choosing a closed positive Weyl chamber $\fa^+$ of $\fa$. 
Let 
$$\Delta=\{\alpha_1, \cdots, \alpha_k\}$$ be the set of simple roots
 $(\fg, \fa^+)$.

As in the introduction, for $g\in G$,
we denote by $\mu(g)\in $ the unique element in $\mathfrak a^+$ such that
\begin{equation*}
g\in K\exp(\mu(g))K.
\end{equation*}

Let $\Ga<G$ be a Zariski dense discrete subgroup.
We denote by $\cal L \subset \fa^+$ the limit cone of $\Ga$, which is the asymptotic cone of $\mu(\Gamma)$:
$$\L=\{\lim t_i \mu(\ga_i)\in \fa^+\text{ for some  $t_i\to 0$ and $\ga_i\in \Ga$}\} .$$

It is a convex cone with non-empty interior \cite{Ben}. The growth indicator function $\psi_\Ga:\fa \to \br\cup\{-\infty\}$ is defined as in \eqref{grow3}. It follows easily from the definition that $\psi_\Ga$ does not depend on 
 the choice of a norm on~$\fa$.

Quint showed the following:
\begin{Thm} \cite[Thm. IV.2.2]{Quint1} \label{growth} The growth indicator function $\psi_\Gamma$ is concave, upper semi-continuous, and satisfies
$$\L= \{u\in \fa^+: \psi_\Gamma(u)>-\infty\}.$$
Moreover, $ {\psi_\Gamma}(u)$ is non-negative on  $\L$ and positive on $\op{int}\L$.
\end{Thm}

\begin{lemma}  \cite[Lem. III.1.3]{Quint1}\label{lem.dom} \label{lem.critless}
Let $F$ be a 
continuous function on $\fa^+$ satisfying $F(tu)=tF(u)$ for all $t\geq 0$ and $u\in\fa$.
If $F(u)>\psi_\Ga(u)$ for all $u\in\fa-\{0\}$, then
$$
\sum_{\ga\in\Ga}e^{-F(\mu(\ga))}<\infty.
$$
Moreover, we have $\delta_{\Ga, F}<1$.
\end{lemma}

\begin{proof} Convergence of the series is shown in \cite[Lem. III.1.3]{Quint1}, and in particular $\delta_{\Ga,F} \le 1$. 
To obtain the strict inequality, we claim that there exists $0< \varepsilon <1$ such that \be \label{ww} (1 - \varepsilon) F > \psi_{\Ga} \quad \mbox{on } \fa - \{0\}.\ee 
Since $\psi_\Ga=-\infty$ outside $\L$ and both
$F$ and $\psi_\Ga$ are homogeneous functions,
it suffices to prove \eqref{ww} on $\{\|v\|=1, v\in \cal L\}$.
Since $ \psi_\Ga\ge 0$ on $\L$, we have $F>0$ on $\L-\{0\}$.
Hence the claim now follows because
$\frac{\psi_\Ga}{F}$ is upper semi-continuous and thus achieves its maximum on any compact set.
\end{proof}

 We denote by $\fa^*$ the set of all linear forms on $\fa$.
\begin{Def} A linear form $\alpha\in \fa^*$ is called {\it tangent} to $\psi_\Ga$ at $u \in \fa - \{0\}$ if $\alpha\ge \psi_\Ga$ and $\alpha(u)=\psi_\Ga(u)$.
\end{Def} 

Consider the following dual cone of the limit cone $\cal L$: \be\label{dual2} \L^\star:=\{\alpha \in \fa^*: \alpha(v) \ge 0 \mbox{ for all } v \in \L \}.\ee 
Observe that the set of all positive roots
is contained in $\L^\star .$

Note that the interior of $\L^\star$ is given as
$$\inte\L^\star=\{\alpha \in \fa^* : \alpha(v) > 0 \mbox{ for all }  v \in \L - \{0\} \}.$$

For any $\alpha \in  \L^{\star}$, 
we set
$$\delta_\alpha=\delta_{\Ga, \alpha}.$$

\begin{lemma} \label{finite}
If $\alpha \in \inte \L^{\star}$,  then 
$$\delta_\alpha\le \sup_{v\in \L-\{0\}}\frac{\psi_\Ga(v)}{\alpha (v)}<\infty. $$

\end{lemma} 
\begin{proof}
Let $\kappa:=\sup_{v\in \L-\{0\}}\frac{\psi_\Ga(v)}{\alpha (v)}$.
Since $\alpha >0$ on $\L-\{0\}$, $0\le \kappa =\sup_{v\in \L-\{0\}}\frac{\psi_\Ga(v)}{\alpha (v)}<\infty $ is well-defined.
Since
$\psi_\Ga< (\kappa +\e) \alpha$ on $\fa - \{0\}$ for any $\e>0$,
 we have, by Lemma \ref{lem.dom}, that $\delta_{(\kappa+\e)\alpha} < 1$.
Hence $\delta_\alpha < \kappa +\e$.
Since $\e>0$ is arbitrary, we get $\delta_\alpha \le \kappa$.
\end{proof}

\begin{Thm} \label{thm.bigger} Let $\Ga<G$ be a Zariski dense discrete subgroup.
For any non-zero
$\alpha \in \L^{\star}$ with $\delta_\alpha<\infty$, 
the linear form
$$T_{\alpha} := \delta_{\alpha}\alpha$$ is tangent to
$\psi_\Ga$ and $\delta_\alpha>0$.
In particular, for any subset $S\subset\inte \L^\star$,
$$\psi_\Ga \le \inf_{\alpha\in S} T_\alpha  .$$
\end{Thm}

\begin{proof}
Fix any norm $\| \cdot \|$ on $\fa$ and we use this norm in the definition of $\psi_{\Ga}$.
We first claim 
\be\label{first} \psi_{\Ga}(v) \le \delta_{\alpha} \alpha(v) \quad \mbox{for all } v \in  \inte \L.\ee  Fix 
$v \in \inte \L$
and $\varepsilon > 0$. We then consider 
$$\cal C_{\varepsilon}(v) = \left\{w \in \fa :  \alpha(w) > 0 \mbox{ and } \left| \frac{\|w\|}{\alpha(w)} - \frac{\|v\|}{\alpha(v)} \right| < \varepsilon  \right\};$$ 
since $\alpha(v)>0$, this is a well-defined 
open cone containing $v$.
Therefore by the definition of $\psi_\Ga$,
we have \be\label{tau} \psi_{\Ga}(v) \le \|v\| \tau_{\cal C_{\varepsilon}(v)}.\ee
Observe that for any $s\ge 0$,
$$\begin{aligned}
\sum_{\ga \in \Ga, \mu(\ga) \in \cal C_{\varepsilon}(v)} e^{-s \| \mu(\ga)\|} & \le \sum_{\ga \in \Ga, \mu(\ga) \in \cal C_{\varepsilon}(v)} e^{-s \alpha(\mu(\ga)) \left( \frac{\|v\|}{\alpha(v)} - \varepsilon \right)} \\
& \le \sum_{\ga \in \Ga}  e^{-s \alpha(\mu(\ga)) \left( \frac{\|v\|}{\alpha(v)} - \varepsilon \right)}.
\end{aligned}$$
Since $\tau_{\cal C_{ \varepsilon}(v)}$ is the abscissa of convergence of the series $$\sum_{\small \ga \in \Ga, \mu(\ga) \in \cal C_{ \varepsilon}(v)} e^{-s \| \mu(\ga) \|} ,$$ 
it follows from the definition of $\delta_\alpha$ that 
$$\tau_{\cal C_{ \varepsilon}(v)} \le \frac{\delta_{\alpha}}{\|v\|{\alpha(v)}^{-1} - \varepsilon} = \frac{\delta_{\alpha}\alpha(v)}{\|v\| - \varepsilon \alpha(v)}.$$
Together with \eqref{tau}, we have $$\psi_{\Ga}(v) \le \|v\| \frac{\delta_{\alpha}\alpha(v)}{\|v\| - \varepsilon \alpha(v)}.$$
Since $\varepsilon > 0$ is arbitrary, we get $$\psi_{\Ga}(v) \le \delta_{\alpha}\alpha(v).$$
This proves the claim \eqref{first}.

We now claim that the inequality \eqref{first} also holds
 for any $v$ in the boundary $\partial \L$.
 Choose any $v_0 \in \inte \L$. From the concavity of $\psi_{\Ga}$, we have $$t\psi_{\Ga}(v_0) + (1-t)\psi_{\Ga}(v) \le \psi_{\Ga}(tv_0 + (1-t)v) \quad \mbox{for all } 0 < t < 1.$$ Since $\L$ is convex, $t v_0 + (1-t) v \in \inte \L$ for all $0 < t < 1$. As we have already shown $\psi_\Ga \le T_\alpha$ on $\inte\L$, we get $$t \psi_{\Ga}(v_0) + (1 - t) \psi_{\Ga}(v) \le T_{\alpha}(tv_0 + (1-t)v) \quad \mbox{for all } 0 < t < 1.$$ By sending $t \to 0^+$, we get $$\psi_{\Ga}(v) \le T_{\alpha}(v). $$ 
Since $\psi_\Ga=-\infty$ outside $\L$, we have established
$\psi_\Ga\le T_\alpha$.  It remains to show that $\psi_\Ga(v)=T_\alpha(v)$ for some  $v \in \fa - \{0\}$. Suppose not, i.e.,
$\psi_\Ga <\delta_\alpha  \alpha $ on $\fa - \{0\}$.
By Lemma \ref{lem.critless}, the abscissa of convergence of the series \be \label{eqn.series}
 \sum_{\ga \in \Ga} e^{-s \delta_{\alpha} \alpha(\mu(\ga))}
\ee is strictly less than $1$. However the abscissa of convergence of the series \eqref{eqn.series} is equal to $1$ by the definition of $\delta_\alpha$. 
 Therefore we have obtained a contradiction.

Note that this implies $\delta_\alpha>0$ since $\psi_\Ga>0$
on $\inte \L$, which is non-empty by Zariski density hypothesis by Theorem \ref{growth}. The last part of the theorem follows from Lemma \ref{finite}.
\end{proof}

\begin{rmk}
We also note the following lower bound for $\psi_\Ga$: 
let $T_\ell \in \L^{\star}$, $\ell\in I$, be a finite collection of linear forms which are tangent to $\psi_\Ga$ at some $v_\ell\in \L-\{0\}$.
Then the concavity property of $\psi_\Ga$ implies that
for any $v=\sum_{\ell\in I} c_\ell v_\ell$ with $c_\ell \ge 0$,
$$\sum_{\ell \in I} c_\ell T_{\ell}(v_\ell) \le \psi_{\Ga}(v).$$
\end{rmk}

\noindent{\bf Proof of Theorem \ref{m}} Note that $\Delta\subset \L^\star$.
Hence this follows from Theorem \ref{thm.bigger} by taking the minimum over all simple roots $\alpha_i \in \Delta$ with $\delta_{\alpha_i}<\infty$.

\medskip

We also note the following corollary of Theorem \ref{thm.bigger}:
\begin{cor} \label{thm.bigger0}  Let $\Ga<G$ be a Zariski dense discrete subgroup.
 For any $\alpha \in \inte \L^{\star}$,  
we have
$$ 0<\delta_\alpha =\max_{v\in \L-\{0\}}\frac{\psi_\Ga(v)}{\alpha (v)} <\infty.$$
\end{cor}
\begin{proof} 
By Lemma \ref{finite}, $\delta_{\alpha} < \infty$. Hence Theorem \ref{thm.bigger} implies $\psi_{\Ga} \le \delta_{\alpha} \alpha$ and $\psi_{\Ga}(v) = \delta_\alpha \alpha(v)$ for some $v\ne 0$. This implies the claim.
\end{proof}

By the following theorem, the above corollary applies to $\alpha \in \theta$
for $\theta$-Anosov subgroups.
\begin{Thm}[{\cite{GGKW}, \cite{KLP1}}] \label{klp}
If $\Ga$ is $\theta$-Anosov, then
$$\theta\subset \inte \L^\star .$$
In particular, if $\Ga$ is $\Delta$-Anosov, then
\be\label{li} \L\subset \inte \fa^+ \cup\{0\}.\ee 
\end{Thm}

\section{Proof of Theorem \ref{m4}} \label{sec.dual}

In this section, let 
$$\text{$\Ga < G$ be a Zariski dense $\Delta$-Anosov subgroup,}$$ as defined in the introduction \eqref{eqn.Anosovdef}.

By  Quint's duality lemma \cite[Lem. 4.3]{Q4} and the works of Quint \cite{Q4}, Sambarino \cite[Lem. 4.8]{Samb3} and Potrie-Sambarino \cite[Prop. 4.6 and 4.11]{SP}, which is based on the work \cite{BCLS},
we have the following fundamental properties of $\Ga$:

\begin{Thm} \label{SPS}   On $\inte \L$,
    $\psi_\Ga$ is analytic, strictly concave, and vertically tangent on $\partial \L$. 
\end{Thm}
The vertical tangency of $\psi_\Ga$ on $\partial \L$ means that there are no linear forms which are tangent to $\psi_\Ga$ at a point of $\partial \L$.
 
In the following, we fix a  norm on $\fa$ induced from an inner product
$\langle \cdot, \cdot \rangle$ which is non-negative on $\fa^+$, i.e., $\langle u, v \rangle \ge 0$ for all $u, v \in \fa^+$. We denote by
$\nabla \psi_\Ga  (u)\in \fa$ the gradient of $\psi_\Ga$ at $u$ so that 
$ d(\psi_\Ga)_u(v)=\langle \nabla \psi_\Ga(u), v\rangle$ for all  $v\in \fa-\{0\}$.

The following theorem was first observed by Quint for Schottky groups \cite{Q4} and is deduced from Theorem \ref{SPS} in general:
\begin{Thm}[{{\cite[Coro. 7.8]{ELO}} \cite[Prop. 4.4]{LO}}] \label{p2} 
Let $u\in \inte \L$.
\begin{enumerate} 
\item 
    There exists a unique $\psi_u\in  \fa^*$ which is tangent to $\psi_\Ga$ at $u$.
    \item We have $\psi_u\in \inte \L^\star$ and
    \be\label{pu} \psi_u(\cdot) = \langle \nabla \psi_{\Ga}(u), \cdot \rangle =d(\psi_\Ga)_{u}.
    \ee 
   \item The map $u\mapsto \psi_u$ induces a bijection between directions in $\inte \L$ and directions in $\inte \L^\star$.
\item We have $\delta_{\psi_u}=1$.

\end{enumerate}

\end{Thm}

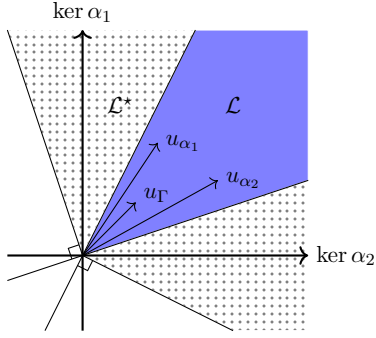
\begin{figure}[ht]
\centering
\begin{tikzpicture}[scale=1, every node/.style={scale=0.8}]

\filldraw[draw= white, pattern=dots, pattern color=gray] (-1, 3) -- (0, 0) -- (2, -1) -- (3, -1) -- (3, 3) -- (-1, 3);
\filldraw[draw= white, fill = blue!50, opacity=0.2] (3/2, 3) -- (0, 0) -- (3, 1) -- (3, 3) -- (3/2, 3);

\draw[->, thick] (-1, 0) -- (3, 0);
\draw[->, thick] (0, -1) -- (0, 3);

%\draw[->] (0, 0) -- (1, 0);
%\draw[->] (0, 0) -- (0, 1);
%\draw (1, 0) node[below] {$e_1$};
%\draw (0.1 ,1.2) node[left] {$e_2$};
\draw (3, 0) node[right] {$\ker \alpha_2$};
\draw (0, 3) node[above] {$\ker \alpha_1$};

\draw[->] (0, 0) -- (1, 1.5);
\draw (1, 1.5) node[right] {$u_{\alpha_1}$};

\draw[->] (0, 0) -- (1.8, 1);
\draw (1.8, 1) node[right] {$u_{\alpha_2}$};

\draw[->] (0, 0) -- (1/1.414, 1/1.414);
\draw (1/1.414, 0.1 + 1/1.414) node[right] {$u_{\Ga}$};

\draw (-1, -1/3) -- (3, 1);
\draw (-1/2, -1) -- (3/2, 3);

\draw (-1, 3) -- (0, 0);
\begin{scope}[scale=0.5]

\draw (-0.3/3.16227766017, 0.9/3.16227766017) -- (-1.2/3.16227766017, 0.6/3.16227766017) -- (-0.9/3.16227766017, -0.3/3.16227766017);

\end{scope}

\draw (0, 0) -- (2, -1);
\begin{scope}[scale=0.5]
\draw (0.6/2.2360679775, -0.3/2.2360679775) -- (0.3/2.2360679775, -0.9/2.2360679775) -- (-0.3/2.2360679775, -0.6/2.2360679775);
\end{scope}

%\draw (3, 1) node[right] {$d_-$};
%\draw (3/2, 3) node[above] {$d_+$};
%\draw (-1, 3) node[above] {$-1/d_-$};
%\draw (2, -1) node[below] {$-1/d_+$};

\draw (2, 2) node {$\L$};
\draw (0.5, 2) node {$\L^{\star}$};

\end{tikzpicture}

\caption{\small{Limit cone and its dual cone.}} \label{dualcone}
\end{figure}

We deduce the following from the above two theorems:
\begin{prop}\label{dual}
Consider the map $\inte \L\to \inte\L^\star$ given by $u\mapsto \alpha_u$ where
    $$ \alpha_u :=\frac{\psi_u}{\psi_{\Ga}(u)} .$$
\
\begin{enumerate}
    \item  The map $u\mapsto \alpha_u$ is a bijection. 
\item Its inverse map $\inte\L^\star \to\inte \L$ is given by $\alpha\to u_\alpha$
where  $u_\alpha\in \inte \L$ is the unique vector such that $\nabla \psi_\Ga(u_\alpha)$ is perpendicular to $\ker \alpha$ and $$ \alpha ( u_\alpha)=~1.$$
We also have
\be\label{pp10}\psi_{\Ga}(u_{\alpha})=\max_{ v\in\cal L, \alpha ( v) =1 } \psi_\Ga(v). \ee 
\end{enumerate}
\end{prop}
\begin{proof} For $t>0$, $\psi_{tu}=\psi_u$ and $\psi_\Ga(tu)=t\psi_\Ga(u)$;
hence $\alpha_{tu}=t^{-1}\alpha_u$. Therefore (1) follows from Theorem \ref{p2}.

Let $\alpha \in \inte \L^{\star}$. 
 Let $u_{\alpha} \in \inte \L$ be the vector given by the relation
 $\alpha_{u_\alpha}=\alpha$, that is,
$\alpha  = \frac{\psi_{u_\alpha}}{ \psi_{\Ga}(u_{\alpha})} $.
By the definition of $\psi_{u_\alpha}$ given in \eqref{pu}, 
$\nabla \psi_{\Ga}(u_{\alpha})$ is perpendicular to $\ker \alpha$, and
$$ \alpha ( u_\alpha) =\frac{\psi_{u_\alpha}(u_\alpha)}{\psi_{\Ga}(u_{\alpha})} = \frac{ \psi_{\Ga}(u_{\alpha}) }{\psi_{\Ga}(u_{\alpha})}=1.$$
 To show the uniqueness, suppose that $v\in \inte\L$
 is a vector such that
 $\nabla \psi_{\Ga}(v)$ is parallel to
 $\nabla \psi_{\Ga}(u_\alpha)$ and $\alpha(v)=1$.
The strict concavity of $\psi_{\Ga}$ on $\inte \L$ as in Theorem \ref{SPS} implies that $v$ must be parallel to $u_\alpha$.  Since
 $\alpha(v)=\alpha(u_\alpha)=1$, it follows that $v=u_\alpha$.

 Observe that for any $v \in \L$ with $\alpha( v)= 1$, we have $$\psi_{\Ga}(v) \le \psi_{u_{\alpha}}(v) =  \psi_{\Ga}(u_{\alpha})\alpha ( v) = \psi_{\Ga}(u_{\alpha})=\psi_\Ga(u_\alpha).$$
Since $\alpha(u_{\alpha}) = 1$, this implies \eqref{pp10}.
\end{proof}

\begin{figure}[ht]
\centering
\begin{tikzpicture}[scale=1, every node/.style={scale=0.7}]

\filldraw[draw= white, fill = blue!50, opacity=0.2] (3/2, 3) -- (0, 0) -- (3, 1) -- (3, 3) -- (3/2, 3);

\draw[->, thick] (-1, 0) -- (3, 0);
\draw[->, thick] (0, -1) -- (0, 3);

\draw (-1, -1/3) -- (3, 1);
\draw (-1/2, -1) -- (3/2, 3);
\draw (-1, 1.5) node[left] {parallel to $\ker \alpha$};
\draw (3, 1.5) node[right] {\color{white} parallel to $\ker \alpha$};

% y = -1/3x + 1 -> (-1, 4/3), (3, 0)

\draw (-1, 4/3) -- (3, 0);
\draw[thick] (3/7, 6/7) -- (3/2, 1/2);

\draw (3/2+0.1, 3) .. controls (1.25+0.1, 2.5) and (0.55, 0.15 + 2/3) .. (1, 2/3) .. controls (1.3, -0.1 + 2/3) and (2.7, 0.9+0.1) .. (3, 1+0.1);

\draw[->] (-1/2, 7/6) -- (-1/6, 13/6);
\draw (-1/2 - 0.15, 7/6 + 0.05) -- (-1/2 - 0.15 + 0.05, 7/6 + 0.05 + 0.15) -- (-1/2 - 0.15 + 0.05 + 0.15, 7/6 + 0.05 + 0.15 - 0.05);

\begin{scope}[shift={(1, 2/3)}]
\begin{scope}[scale=0.5]
\begin{scope}[shift={(1/2, -7/6)}]
    \draw[->] (-1/2, 7/6) -- (-1/6, 13/6);
    \draw (-1/2 - 0.15, 7/6 + 0.05) -- (-1/2 - 0.15 + 0.05, 7/6 + 0.05 + 0.15) -- (-1/2 - 0.15 + 0.05 + 0.15, 7/6 + 0.05 + 0.15 - 0.05);
\end{scope}
\end{scope}
\end{scope}

\draw[->] (0, 0) -- (0.97, 1.94/3);
\draw (0.5, 0.5) node {$u_{\alpha}$};

\draw (1.2, 1.2) node[right] {$\nabla \psi_{\Ga}(u_{\alpha})$};

\draw (2, 3) node[above] {$\{\psi_{\Ga} = \psi_{\Ga}(u_{\alpha})\}$};
%\draw (-1/6, 13/6) node[left] {$\alpha$};

\end{tikzpicture}

\caption{\small{From $\alpha$ to $u_{\alpha}$}} \label{grad}
\end{figure}
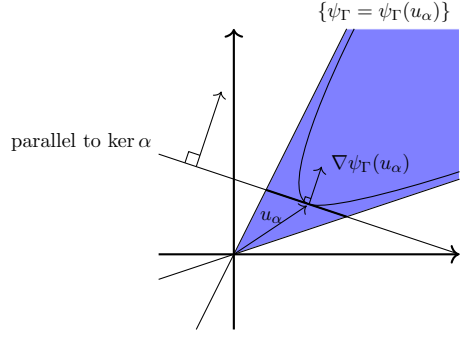

\begin{Thm} \label{alpha}\label{tang}\label{tent}
    For any $\alpha\in \inte\L^\star$,
    we have $$\delta_{\alpha} =\psi_{\Ga}(u_{\alpha}) \quad \text{ and } \quad \psi_{u_\alpha}=\delta_\alpha \alpha .$$
\end{Thm}
\begin{proof}
The first claim follows from \eqref{pp10} and Corollary \ref{thm.bigger0}. Since $\psi_{u_{\alpha}} = \psi_{\Ga}(u_{\alpha}) \alpha$ by Proposition \ref{dual}, the first claim implies the second. 
\end{proof}

 \subsection*{Proof of Theorem \ref{thm.Anosov}}
For \eqref{(2)}, we claim that $v_i:=u_{\alpha_i}$ satisfies the claim. 
 By Proposition \ref{dual}, we have $u_{\alpha_i} \in \inte \L$ and it satisfies $\alpha_i(u_{\alpha_i}) = 1$. By Lemma \ref{alpha}, $\psi_{\Ga}(u_{\alpha_i}) = \delta_{ \alpha_i}$. 
 The uniqueness follows easily from the strict concavity of $\psi_\Ga$
 (Theorem \ref{SPS}).

For \eqref{(3)}, suppose that
for some $v \in \fa$ and $1 \le i \le k$,
we have $\psi_{\Ga}(v) = \delta_{\alpha_i} \alpha_i(v)$. Since $\psi_{u_{\alpha_i}}=\delta_{\alpha_i} \alpha_i$ is a tangent form to $\psi_\Ga$ at $u_{\alpha_i}$,
it follows again from the strict concavity of $\psi_\Ga$
and the vertical tangency property (Theorem \ref{SPS} that $v$ is parallel to  $u_{\alpha_i}$.

By Theorem \ref{m}, we have \be \label{eqn.pi2} \psi_{\Ga}(v) \le \min_{1 \le i \le k} \delta_{ \alpha_i} \alpha_i(v) \le \frac{1}{k} \sum_{1 \le i \le k} \delta_{ \alpha_i} \alpha_i(v).\ee  Suppose that $\psi_{\Ga}(v) = \frac{1}{k} \sum_{1 \le i \le k} \delta_{ \alpha_i} \alpha_i(v)$ for some $v \neq 0$. It then follows from \eqref{eqn.pi2} that
$$ \psi_{\Ga}(v)=\min_{1 \le i \le k} \delta_{ \alpha_i} \alpha_i(v)  = \frac{1}{k} \sum_{1 \le i \le k} \delta_{ \alpha_i} \alpha_i(v).$$ It implies that for all $1 \le i \le k$,
  $$\psi_{\Ga}(v) =   \delta_{ \alpha_i} \alpha_i(v).$$  
Then, as we just have seen, this implies that $v$ is parallel to all
 $ u_{\alpha_i}$,  $1 \le i \le k$.
When $k\ge 2$, this contradicts  Theorem \ref{p2}.
This proves \eqref{(4)}.

\subsection*{Proof of Corollary \ref{innerprod}}
For simplicity, we omit $\|\cdot\|$ in the subscript in this proof, e.g.,
$u_\Ga=u_{\Ga, \|\cdot\|} $. Recall that $\delta_\Ga=\psi_\Ga(u_\Ga)$.
Since $\delta_\Ga=\max_{\|v\|=1}\psi_\Ga(v)$,  \cite[Lem. III.3.4]{Quint1}, applied to $\psi_\Ga$, implies that there exists a tangent form, $\psi_{u_\Ga},
$ to $\psi_{\Ga}$ at $u_{\Ga}$. By the vertical tangent condition
in Theorem \ref{SPS}, it follows that  $u_{\Ga} \in \inte \L$.
Moreover,
 we have $\nabla \psi_{\Ga}(u_\Ga)\in \br_{>0} u_\Ga$ \cite[Lem. 2.24]{ELO}.
Therefore, by Theorem \ref{p2}(2),
there exists $c_0>0$ such that
\be\label{psiu} \psi_{u_\Ga}(\cdot)  =\langle c_0u_\Ga, \cdot \rangle .\ee 

We now claim that 
$$\psi_\Ga (u_\Ga)<T_\Ga (u_\Ga).$$

Suppose not. Then, by Theorem \ref{thm.Anosov}, there exist  $c>0$ and $1 \le i \le k$ such that
 $u_\Ga=c u_{\alpha_i}$ and hence $\psi_{u_\Ga}=\psi_{u_{\alpha_i}} =\delta_{\alpha_i}\alpha_i$. By \eqref{psiu}, it follows that
$\alpha_i(\cdot)  =\langle c_1u_\Ga, \cdot \rangle  $ for some $c_1>0$.

Since $u_\Ga\in \inte \L$, and $\langle, \rangle$ is non-negative on $\fa^+$ by the hypothesis,
the linear form $\langle c_1u_\Ga, \cdot \rangle$ is positive on 
$\fa^+-\{0\}$.
On the other hand, the simple root 
$\alpha_i$ is zero on a wall of $\fa^+$. Therefore we obtained a contradiction.
This finishes the proof.

\medskip

We note that in the above proof, the hypothesis that the norm
$\|\cdot\|$ is induced from 
an inner product was used to deduce that $\psi_{u_{\Ga, \|\cdot\|}}$ is strictly positive on $\fa^+-\{0\}$.

\begin{Rmk} \label{comparePS}
We explain how Theorem \ref{thm.Anosov} can be compared with \cite[Coro. 1.4]{SP}. 
Let $(\fa^+)^{\star} = \{ \alpha \in \fa^* : \alpha(v) \ge 0 \mbox{ for all } v \in \fa^+ \}$ so that $$\inte (\fa^+)^{\star} = \{ \alpha \in \fa^* : \alpha(v) > 0 \mbox{ for all } v \in \fa^+ -\{0\} \}.$$

Recall that \cite[Coro. 1.4]{SP} concerns the Hitchin representations, but their argument applied to our Zariski dense Anosov subgroups yields the following: 
For any $\alpha \in \inte (\fa^+)^\star$,
the quantity $\delta_\alpha$ 
satisfies \be \label{pscoro14}
\delta_{\alpha}\le \frac{1}{\sum_{i = 1}^k a_{i}}
\ee 
where $\alpha = \sum_{i = 1}^k (a_{i} \delta_{\alpha_i}) \alpha_i$; the hypothesis $\alpha \in  \inte (\fa^+)^\star$ is equivalent to $\alpha \neq 0$ and $a_i > 0$ for all $1 \le i \le k$.

On the other hand,
our Theorem \ref{thm.Anosov} says that for all $\alpha \in \inte \L^\star$,
\be \label{ttt} 
\delta_\alpha =\psi_{\Ga}(u_{\alpha}) \le \min_{1 \le i \le k} \delta_{\alpha_i}\alpha_i ( u_{\alpha});\ee
this is equivalent to saying that
for all $v\in \fa$, $\psi_\Ga(v)\le \min_{1 \le i \le k} \delta_{ \alpha_i}\alpha_i ( v)$.

Since \be \label{minineq}
1 =  \alpha ( u_{\alpha}) = \sum_{i = 1}^k a_i  \delta_{\alpha_i}\alpha_i ( u_{\alpha}) \ge \left( \sum_{i = 1}^k a_i  \right) \min_{1 \le i \le k}  \delta_{\alpha_i}\alpha_i( u_{\alpha} )\ee
 we  have $$\min_{1 \le i \le k}  \delta_{\alpha_i}\alpha_i (u_{\alpha} )\le \frac{1}{\sum_{i = 1}^k a_i}$$
where the equality is strict except for one direction of
$u_{\alpha}$ satisfying 
$$ \delta_{\alpha_i} \alpha_i( u_{\alpha}) = \delta_{\alpha_j} \alpha_j( u_{\alpha})  \quad \mbox{for all } i, j = 1, \cdots, k.$$ Therefore our bound \eqref{ttt} is sharper than
the bound \eqref{pscoro14} in addition to the point that it applies
to the optimal cone $\inte \L^\star$, while
\cite[Coro. 1.4]{SP} applies
only for $\alpha \in \inte (\fa^+)^\star$, which is strictly smaller than $\inte \L^\star$.

Both approaches are based on the observation that
the linear forms $\delta_{\alpha} \alpha$'s are tangent to $\psi_\Ga$ for $\alpha \in \Delta$, but 
\cite[Coro. 1.4]{SP} considers these tangent forms as points on the
boundary of the subset  $\cal D=\{\varphi \in \inte (\fa^+)^\star: \delta_\varphi \le 1\}$
and deduce \eqref{pscoro14} from the convexity of $\cal D$, whereas we think of the tangent forms
as functions on $\fa$ and obtain a stronger bound of \eqref{ttt} simply by taking {\it minimum} of these tangent forms over $\alpha\in \Delta$.
\end{Rmk}

\subsection*{Alternate proof of Theorem \ref{thm.Anosov}\eqref{(3)}}
For Anosov subgroups, we present an alternate proof~of 
\be \label{eqn.alt}
\psi_{\Ga} \le T_{\Ga}
\ee for $k \le 3$, using the following ``strip theorem":

\begin{Thm}[Strip theorem] \cite[Thm. 6.3]{BLLO} \label{thm.strip} 
 Let $\Ga$ be a Zariski dense $\Delta$-Anosov subgroup of $ G$.
Let $k=\#\Delta \le 3$ and $v\in \inte \L$. For all sufficiently large $R>0$,
  the abscissa of convergence of the series
 $$\sum_{\gamma\in \Ga, \|\mu(\gamma)-\br v\|\le R} e^{-s \psi_v(\mu(\ga))}$$ is equal to $1$.
\end{Thm}

To show the inequality \eqref{eqn.alt}, we fix $v \in \inte \L$ and $1 \le i \le k$. For $R > 0$, we write $S_R := \{g \in G : \| \mu(g) - \R v \| < R\}$. By Theorem \ref{thm.strip}, there exists $R > 0$ such that the series $\cal D_{R}(s) = \sum_{\ga \in \Ga \cap S_R} e^{-s \psi_v(\mu(\ga))}$ has the abscissa of convergence~$1$. Recalling that $\alpha_i > 0$ on $\inte \L_{\Ga}$, there exists $C > 0$ so that for any $\ga \in S_R$, we have $$\left\| \mu(\ga) - {\alpha_i(\mu(\ga)) \over \alpha_i(v)}v \right\| \le C.$$ It then follows that $$\cal D_R(s) = \sum_{\ga \in \Ga \cap S_R} e^{-s \psi_v(\mu(\ga))} \ll \sum_{\ga \in \Ga \cap S_R} e^{-s {\alpha_i(\mu(\ga)) \over \alpha_i(v)} \psi_{\Ga}(v)} \le \sum_{\ga \in \Ga} e^{-s {\alpha_i(\mu(\ga)) \over \alpha_i(v)} \psi_{\Ga}(v)}.$$ Since the series $\sum_{\ga \in \Ga} e^{-s {\alpha_i(\mu(\ga)) \over \alpha_i(v)} \psi_{\Ga}(v)} $ is finite whenever $s > {\alpha_i(v) \over \psi_{\Ga}(v)} \delta_{\alpha_i}$, we have $1 \le {\alpha_i(v) \over \psi_{\Ga}(v)} \delta_{\alpha_i}$. Hence $$\psi_{\Ga}(v) \le \delta_{\alpha_i} \alpha_i(v).$$ Since  $v \in \inte \L$ and $1 \le i \le k$ are arbitrary, we get $$\psi_{\Ga} \le T_{\Ga} \quad \mbox{on } \inte \L.$$
By the concavity of $\psi_{\Ga}$, this implies $\psi_{\Ga} \le T_{\Ga}$ on $\L$ as well (see the proof of Theorem \ref{thm.bigger}). Since $\psi_{\Ga} = -\infty$ outside $\L$, \eqref{eqn.alt} follows.

\section{Applications to self-joinings} \label{sec.self}
We consider the case when $G = \prod_{i = 1}^k \SO^{\circ}(n_i, 1)$, $n_i \ge 2$, and $\rho_i : \Sigma \to \SO^{\circ}(n_i, 1)$ is a faithful convex cocompact representation with Zariski dense image.  We let
$\Gamma_{\rho}<G$ be the subgroup defined as in \eqref{grho}.
The hypothesis on $\rho_i$'s implies that $\Ga_\rho$ is $\Delta$-Anosov.
We assume $$\text{$k\ge 2$ and $\Ga_\rho$ is Zariski dense in $G$}$$
in the entire section.

\subsection*{Proof of Corollaries \ref{pregap} and \ref{rigid}}
Corollary \ref{pregap} follows since $\delta_{ \alpha_i}=\delta_{\rho_i}=\dim \La_{\rho_i}$.
 For Corollary \ref{rigid}, note that we have $$
\begin{aligned} \delta_{\Gamma_{\rho}, \|\cdot\|_{\op{Euc}}} & < \frac{1}{k} \left( \sum_{i = 1}^k (\dim \La_{\rho_i})^2 \right)^{1/2} \\
& \le \frac{1}{k} \left( k \max_{1 \le i \le k} (\dim \La_{\rho_i})^2 \right)^{1/2}\\& = \frac{1}{\sqrt{k}} \max_{1 \le i \le k} \dim \La_{\rho_i}.
\end{aligned}$$ 
On the other hand, we showed in \cite{KMO}, 
$$\dim \La_{\rho} = \max_i \La_{\rho_i} .$$
Hence $$\delta_{\Gamma_{\rho}, \|\cdot\|_{\op{Euc}}} < \frac{1}{\sqrt{k}} \dim \La_{\rho} .$$

\subsection*{Critical exponent with respect to the $L^1$-metric} Set $\delta_{L^1}:=\delta_{\sum_{i =1 }^k \alpha_i}$, which is the critical exponent of $\Ga_{\rho}$ for the $L^1$-metric $\sum_{i = 1}^k d_i$ on $X = \prod_{i =1 }^k \bH^{n_i}$. We deduce the following from Corollary \ref{pregap}, whose special case when $k=2$ and $\dim \La_{\rho_i}=1$ was proved by Bishop and Steger \cite{BS}:

\begin{cor}  We have
\be\label{one} \delta_{L^1} < \frac{\dim \La_{\rho}}{k} .\ee 
\end{cor}
\begin{proof}
Noting $\alpha:=\sum_{i = 1}^k \alpha_i \in \inte\L^\star$, write $u_{\alpha} =(u_1, \cdots, u_k)\in \inte\L$. Lemma \ref{alpha} and Corollary \ref{pregap} imply
$$\delta_{L^1} =\psi_\Ga (u_{\alpha}) < \frac{1}{k} \sum_{i = 1}^k \dim \La_{\rho_i} u_i \le \frac{\max_i \dim \La_{\rho_i}}{k} \sum_{i = 1}^k u_i.$$
Since $\alpha(u_\alpha)=\sum_{1 \le i \le k} u_i = 1$ by Lemma \ref{dual}(2) and $\max_i \dim \La_{\rho_i} = \dim \La_{\rho}$ by \cite{KMO}, we get the desired inequality.
\end{proof}

\subsection*{Geodesic stretching between two hyperbolic manifolds}
When $k = 2$,  the limit cone $\L$ of $\Ga_{\rho}$ can also be described as
$$\cal L:=\{(v_1,v_2)\in \br_{\ge 0}^2: d_-v_1\le v_2\le d_+v_1\} $$
where $d_+$ and $d_-$
 are respectively the maximal and minimal geodesic stretching constants of $\rho_2$ relative to $\rho_1$:
$$ d_+(\rho_1,\rho_2) =\sup_{\sigma\in \Sigma-\{e\} }\frac {\ell_2(\s)}{\ell_1(\s)} \quad \mbox{and} \quad d_-(\rho_1,\rho_2) =\inf_{\s\in \Sigma-\{e\}}\frac{\ell_2(\s)}{\ell_1(\s)}$$ 
where
$\ell_i(\s)$ denotes the length of the closed geodesic in the hyperbolic manifold $\rho_i(\Delta)\ba \bH^{n_i}$ corresponding to $\rho_i(\s)$ (cf. \cite{Bu}, \cite{Ben}).

Thurston \cite{Th} showed that the maximal geodesic stretching constant is always strictly bigger than $1$ for finite-area hyperbolic surfaces. (See also \cite{GK}). Theorem \ref{m4}  implies the following corollary; this was
already observed by Burger \cite[Thm. 1 and its Coro.]{Bu} and generalizes a theorem of Thurston \cite[Thm. 3.1]{Th}:

\begin{cor}\label{th} 
We have
$$d_-(\rho_1, \rho_2) < \frac{\dim \La_{\rho_1}}{\dim \La_{\rho_2}} < d_+(\rho_1, \rho_2).$$
\end{cor}
\begin{proof}
By Theorem \ref{m4}, \be \label{observe}
\psi_{\Ga} \le \min (\delta_{ 1} \alpha_1, \delta_{ 2} \alpha_2).
\ee 
 By Theorem \ref{tent}, we have $\psi_{\Ga}(u_{\alpha_1}) = \delta_{1} \alpha_1(u_{\alpha_1})$. Hence
 $$\delta_{1} \alpha_1(u_{\alpha_1}) \le \min (\delta_{1}\alpha_1(u_{\alpha_1}), \delta_{ 2} \alpha_2(u_{\alpha_1})),$$
 which implies $\delta_{{1}} \alpha_1(u_{\alpha_1}) \le \delta_{2} \alpha_2(u_{\alpha_1})$. Therefore, $$\frac{\delta_{1}}{\delta_{ 2}} \le \frac{\alpha_2(u_{\alpha_1})}{\alpha_1(u_{\alpha_1})}.$$
Similarly, we have $\delta_{ 2} \alpha_2(u_{\alpha_2}) \le \min(\delta_{ 1} \alpha_1(u_{\alpha_2}), \delta_{2} \alpha_2(u_{\alpha_2}))$, and hence $$\frac{\alpha_2(u_{\alpha_2})}{\alpha_1(u_{\alpha_2})} \le \frac{\delta_{1}}{\delta_{ 2}}.$$ Since $\dim \La_{\rho_i} = \delta_{ i}$ for $i = 1, 2$ by Patterson \cite{Pa} and Sullivan \cite{Su}, we now have $$\frac{\alpha_2(u_{\alpha_2})}{\alpha_1(u_{\alpha_2})} \le \frac{ \dim \La_{\rho_1}}{\dim \La_{\rho_2}} \le \frac{ \alpha_2(u_{\alpha_1})}{\alpha_1(u_{\alpha_1})}.$$ Since $u_{\alpha_1}, u_{\alpha_2} \in \inte \L$, $d_-(\rho_1, \rho_2) < \frac{ \alpha_2(u_{\alpha_2})}{\alpha_1(u_{\alpha_2})}$ and $\frac{\alpha_2(u_{\alpha_1})}{\alpha_1(u_{\alpha_1})} < d_+(\rho_1, \rho_2)$.
It completes the proof.
\end{proof}

\end{document}